\documentclass[a4,draft,fleqn]{article}
\usepackage{fullpage}
\usepackage{amsmath}
\usepackage{amssymb}
\begin{document}

\title{ Normalization of Hamiltonian in the Generalized Photogravitaional Restricted Three Body Problem with Poynting-Robertson
Drag}
\author{B.S.Kushvah$^1$,  J.P. Sharma$^2$ and  B.Ishwar$^3$ \\\\
1. JRF  DST Project, 2.  Co-P.I   DST Project  and  3. P.I. DST Project\\
    University Department of Mathematics,\\
    B.R.A. Bihar University Muzaffarpur-842001,
    Bihar (India)\\Email:ishwar\_ bhola@hotmail.com, bskush@hotmail.com}
\maketitle 
\noindent 
\begin{abstract}
 We have performed normalization of Hamiltonian in the generalized photogravitational restricted three body problem with Poynting-Robertson drag. In this problem we have taken bigger primary as source of radiation and smaller primary as an oblate spheroid. Wittaker method is used to transform
the second order part of the Hamiltonian into the normal form. 

\end{abstract}
\section{Introduction}

The restricted three body problem describes the motion of an
infinitesimal mass moving under the gravitational effect of the two
finite masses, called primaries, which move in circular orbits
around their centre of mass on account of their mutual attraction
and the infinitesimal mass not influencing the motion of the
primaries. The classical restricted three body problem is
generalized to include the force of radiation pressure, the
Poynting-Robertson effect and oblateness effect.

J. H. Poynting (1903) considered the effect of the absorption and
subsequent re-emission of sunlight by small isolated particles in
the solar system. His work was later modified by H. P. Robertson
(1937) who used a precise relativistic treatments of the first order
in the ratio of he velocity of the particle to that of light.

The effect of radiation pressure and P-R. drag in the restricted
three body problem has been studied by Colombo {\it et al.} (1966),
Chernikov Yu. A. (1970) and Schuerman (1980) who discussed the
position as well as the stability of the Lagrangian equilibrium
points when radiation pressure, P-R drag force are included. Murray
C. D. (1994) systematically discussed the dynamical effect of
general drag in the planar circular restricted three body problem,
Liou J. C. {\it et al.} (1995) examined the effect of radiation
pressure, P-R drag and solar wind drag in the restricted three body problem.

Moser's conditions (1962), Arnold's theorem (1961) and Liapunov's
theorem (1956) played a significant role in deciding the nonlinear
stability of an equilibrium point. Applying Arnold's theorem (1961),
Leontovic (1962) examined the nonlinear stability of triangular
points. Moser gave some modifications in Arnold's theorem. Then
Deprit and Deprit (1967) investigated the nonlinear stability of
triangular points by applying Moser's modified version of Arnold's
theorem (1961).

Bhatnagar and Hallan (1983) studied the effect of perturbations on
the nonlinear stability of triangular points. Maciejewski and
Gozdziewski (1991) described the normalization algorithms of
Hamiltonian near an equilibrium point. Niedzielska (1994) investigated the nonlinear
stability of the libration points in the photogravitational
restricted three body problem. Mishra P. and Ishwar B.(1995) studied second order normalization in the generalized restricted problem of three bodies, smaller primary being an oblate spheroid. Ishwar B.(1997) studied nonlinear stability in the generalized restricted three body problem.

In this paper  normalization of Hamiltonian is performed   in the
generalized photogravitaional restricted three body problem with
Poynting-Robertson drag.  Wittaker method is used to transform
the second order part of the Hamiltonian into the normal form.

\section{Location of Triangular Equilibrium Points}
\newcommand{\abc}{{\Bigl(1+\frac{5}{2}A_2\Bigr) }}
\newcommand{\aac}{{(1-\frac{A_2}{2}) }}
\newcommand{\amc}{{(1-\mu) }}
\newcommand{\adc}{{\frac{\delta ^2}{2} }}
\newcommand{\zab}{{\Bigl(1+\frac{5A_2}{2r^2_2}\Bigr)}}
\newcommand{\zw}{{\frac{W_1}{r^2_1}}}
\newcommand{\zwf}{{\frac{W_1}{r^4_1}}}
\newcommand{\zx}{{(x+\mu)}}
\newcommand{\zox}{{(x+\mu-1)}}
\newcommand{\zd}{{\displaystyle}}
\newcommand{\zabs}{{\Bigl(1+\frac{5A_2}{2{r^2_2}_*}\Bigr)}}
\newcommand{\zws}{{\frac{W_1}{{r^2_1}_*}}}
\newcommand{\zwfs}{{\frac{W_1}{{r^4_1}_*}}}
\newcommand{\zxs}{{(x_*+\mu)}}
\newcommand{\zoxs}{{(x_*+\mu-1)}}
\newcommand{\az}{{1,0}}
\newcommand{\za}{{0,1}}
\newcommand{\zae}{A_2\epsilon}
\newcommand{\zwe}{nW_1\epsilon}
Equations of motions are
\begin{align}
\ddot{x}-2n\dot{y}&=U_x ,\quad\text{where},\quad U_x=\frac{\partial{U_1}}{\partial{x}}-\frac{W_{1}N_1}{r^2_1}\\
\ddot{y}+2n\dot{x}&=U_y,\hspace{.85in}U_y=\frac{\partial{U_1}}{\partial{y}}-\frac{W_{1}N_2}{r^2_1}\\
U_1&=\zd{\frac{n^2(x^2+y^2)}{2}}+\frac{\amc{q_1}}{r_1}+\frac{\mu}{r_2}+\frac{\mu{A_2}}{2r^3_2}
\end{align}
\begin{gather*}
r^2_1=\zx^2+y^2,\quad  r^2_2=\zox^2+y^2,\quad n^2=1+\frac{3}{2}A_2,\\
N_1=\frac{\zx[\zx\dot{x}+y\dot{y}]}{r^2_1}+\dot{x}-ny,\quad
N_2=\frac{y[\zx\dot{x}+y\dot{y}]}{r^2_1}+\dot{y}+n\zx
\end{gather*}
$W_1=\frac{(1-\mu)(1-q_1)}{c_d}$,
$\mu=\frac{m_2}{m_1+m_2}\leq\frac{1}{2}$, $m_1,m_2$ be  the  masses of
the primaries, $A_2=\frac{r^2_e-r^2_p}{5r^2}$ be the oblateness
coefficient, $r_e$ and$r_p$ be the equatorial and polar radii
respectively $r$ be the distance between primaries,
$q=\bigl(1-\frac{F_p}{F_g}\bigr)$ be the mass reduction factor
expressed in terms of the particle's radius $a$, density $\rho$ and
radiation pressure efficiency factor $\chi$ (in the C.G.S.system)
i.e., $q=1-\zd{\frac{5.6\times{10^{-5}}\chi}{a\rho}}$. Assumption
$q=constant$ is equivalent to neglecting fluctuation in the beam of
solar radiation and the effect of solar radiation, the effect of the
planet's shadow, obviously $q\leq1$. Triangular equilibrium points
are given by $U_x=0,U_y=0,z=0,y\neq{0}$, then we have
\begin{align}
 x_*&=x_0\Biggl\{1-\zd{\frac{nW_1\bigl[\amc\abc+\mu\aac\adc\bigr]}{3\mu\amc{y_0 x_0}}}-\adc\frac{A_2}{x_0}\Biggr\} \label{eq:1x}\\
y_*&=y_0{\Biggl\{1-\zd{\frac{nW_1\delta^2\bigl[2\mu-1-\mu(1-\frac{3A_2}{2})\adc+7\amc\frac{A_2}{2}\bigr]}{3\mu\amc{y^3_0}}}-\zd{\frac{\delta^2\bigl(1-\adc)A_2}{y^2_0}}\Biggr\}^{1/2}
}
\end{align}
where
$x_0=\adc-\mu$, $y_0=\pm\delta\bigl(1-\frac{\delta^2}{4}\bigr)^{1/2}$
and $\delta=q^{1/3}_1$, as in preprint Kushvah \& Ishwar (2006)

\section{ Normalization of Hamiltonian}
We used  Whittakar (1965) mathod  for the transformation of $H_2$
into normal form  The Lagrangian function of the problem can be
written as
\begin{align}
L&=\frac{1}{2}(\dot{x}^2+\dot{y}^2)+n(x\dot{y}-\dot{x}y)+\frac{n^2}{2}(x^2+y^2)+\frac{\amc{q_1}}{r_1}+\frac{\mu}{r_2}+\frac{\mu{A_2}}{2r^3_2}\\\notag
&+W_1\Bigl\{\frac{\zx\dot{x}+y\dot{y}}{2r^2_1}-n
\arctan{\frac{y}{\zx}}\Bigr\}\\\notag
\end{align}
and the Hamiltonian  is $H=-L+P_x\dot{x}+P_y\dot{y}$, where
$P_x,P_y$ are the momenta coordinates given by \[
P_x=\frac{\partial{L}}{\partial{\dot{x}}}=\dot{x}-ny+\frac{W_1}{2r_1^2}\zx,
\quad
P_y=\frac{\partial{L}}{\partial{\dot{y}}}=\dot{y}+nx+\frac{W_1}{2r_1^2}y
\]
For simplicity we suppose  $q_1=1-\epsilon$, with $|\epsilon|<<1$
then coordinates of triangular equilibrium points  can be written in
the form
\begin{align}
x&=\frac{\gamma}{2}-\frac{\epsilon}{3}-\frac{A_2}{2}+\frac{A_2
\epsilon}{3}-\frac{(9+\gamma)}{6\sqrt{3}}nW_1-\frac{4\gamma
\epsilon}{27\sqrt{3}}nW_1 \\
y&=\frac{\sqrt{3}}{2}\Bigl\{1-\frac{2\epsilon}{9}-\frac{A_2}{3}-\frac{2A_2
\epsilon}{9}+\frac{(1+\gamma)}{9\sqrt{3}}nW_1-\frac{4\gamma
\epsilon}{27\sqrt{3}}nW_1\Bigr\}
\end{align}
where $\gamma=1-2\mu$.
 We shift the origin to $L_4$. For that, we change
$x\rightarrow {x_*}+x$ and  $y\rightarrow{y_*}+y$. Let $a=x_*+\mu,
b=y_*$ so that
\begin{align}
a&= \frac{1}{2} \Biggl\{-\frac{2\epsilon}{3}-A_2+\frac{2A_2
\epsilon}{3}-\frac{(9+\gamma)}{3\sqrt{3}}nW_1-\frac{8\gamma
\epsilon}{27\sqrt{3}}nW_1 \bigr\}\\
b&=\frac{\sqrt{3}}{2}\Bigl\{1-\frac{2\epsilon}{9}-\frac{A_2}{3}-\frac{2A_2
\epsilon}{9}+\frac{(1+\gamma)}{9\sqrt{3}}nW_1-\frac{4\gamma
\epsilon}{27\sqrt{3}}nW_1\Bigr\}
\end{align}
Expanding $L$ in power series of $x $ and $y$, we get
\begin{eqnarray}
 L&=&L_0+L_1+L_2+L_3+\cdots \\
H&=&H_0+H_1+H_2+H_3+\cdots =-L+P_x{\dot{x}}+P_y{\dot{y}}
  \end{eqnarray}
  where $L_0,L_1,L_2,L_3 \ldots$ are
\begin{eqnarray}
L_0&=&\frac{3}{2}-\frac{2\epsilon}{3}-\frac{\gamma
\epsilon}{3}+\frac{ 3 \gamma A_2}{4}-\frac{3 A_2 \epsilon}{2}-\gamma
A_2 \notag \\
&&-\frac{\sqrt{3}nW_1}{4}+\frac{2\gamma}{3\sqrt{3}}nW_1-\frac{ n
\epsilon W_1}{3\sqrt{3}}-\frac{23\epsilon n W_1}{54\sqrt{3}}-n
\arctan{\frac{b}{a}}
\end{eqnarray}

\begin{align}
L_1&=\dot{x}\bigl\{-\frac{\sqrt{3}}{2}-\frac{5 A_2
}{8\sqrt{3}}+\frac{7\epsilon A_2}{12\sqrt{3}}+\frac{4
nW_1}{9}-\frac{1
 \gamma nW_1}{18}\bigr\}\notag \\&+ \dot{y}\bigl\{\frac{1}{2}-\frac{\epsilon}{3}-\frac{A_2
}{8}+\frac{\epsilon A_2}{12\sqrt{3}}-\frac{ nW_1}{6\sqrt{3}}+\frac{2
n \epsilon nW_1}{3\sqrt{3}}\bigr\} \notag \\
& -x \bigr\{-\frac{1}{2}+\frac{\gamma}{2}+\frac{9
A_2}{8}+\frac{15\gamma A_2}{8}-\frac{35\epsilon
A_2}{12}-\frac{29\gamma \epsilon  A_2}{12}\notag\\&+
\frac{3\sqrt{3}nW_1}{8}-\frac{2\gamma}{3\sqrt{3}}nW_1-\frac{5
\epsilon n W_1}{12\sqrt{3}}-\frac{7 \gamma \epsilon n
W_1}{4\sqrt{3}}\bigr\}\notag\\&-y
\bigr\{\frac{15\sqrt{3}A_2}{2}+\frac{9\sqrt{3}\gamma
A_2}{8}-2\sqrt{3} \epsilon A_2-2\sqrt{3}\gamma \epsilon A_2-
\frac{nW_1}{8}+\gamma nW_1-\frac{43 \epsilon }{36}nW_1 \bigr\}
\end{align}
\begin{eqnarray}
L_2&=&\frac{(\dot x^2+ \dot y^2)}{2}+n(x\dot y-\dot x y)+
\frac{n^2}{2}(x^2+y^2)-Ex^2-Fy^2-G xy
\end{eqnarray}
\begin{equation}
L_3=-\frac{1}{3!}\left\{x^3T_1+3x^2yT_2+3xy^2T_3+y^3T_4+6T_5\right\}
  \end{equation}
  where
\begin{eqnarray}
E&=&\frac{1}{16}\Bigl\{ 2-6\epsilon- 3A_2-
\frac{31A_2\epsilon}{2}-\frac{(69+\gamma)}{6\sqrt{3}}nW_1+\frac{2(307+75\gamma)
\epsilon}{27\sqrt{3}}nW_1 \notag \\&+&\gamma \bigl\{2\epsilon+12A_2+
\frac{A_2\epsilon}{3}+\frac{(199+17\gamma)}{6\sqrt{3}}nW_1-\frac{2(226+99\gamma)
\epsilon}{27\sqrt{3}}nW_1\bigr\}\Bigr\}
  \end{eqnarray}
  \begin{eqnarray}
F&=&\frac{-1}{16}\Bigl\{ 10-2\epsilon+21A_2-
\frac{717A_2\epsilon}{18}-\frac{(67+19\gamma)}{6\sqrt{3}}nW_1+\frac{2(413-39\gamma)
\epsilon}{27\sqrt{3}}nW_1 \notag \\&+&\gamma \bigl\{6\epsilon-
\frac{293A_2\epsilon}{18}+\frac{(187+27\gamma)}{6\sqrt{3}}nW_1-\frac{4(247+3\gamma)
\epsilon}{27\sqrt{3}}nW_1\bigr\}\Bigr\}
  \end{eqnarray}
   \begin{eqnarray}
G&=&\frac{\sqrt{3}}{8}\Bigl\{ 2\epsilon+6A_2-
\frac{37A_2\epsilon}{2}-\frac{(13+\gamma)}{2\sqrt{3}}nW_1+\frac{2(79-7\gamma)
\epsilon}{27\sqrt{3}}nW_1 \notag \\&-&\gamma
\bigl\{6\epsilon-\frac{\epsilon}{3}+13A_2-
\frac{33A_2\epsilon}{2}+\frac{(11-\gamma)}{2\sqrt{3}}nW_1-\frac{(186-\gamma)
\epsilon}{9\sqrt{3}}nW_1\bigr\}\Bigr\}
  \end{eqnarray}
 \begin{eqnarray} T_1&=&\frac{3}{16}\biggl[\frac{16}{3}\epsilon+6A_2-\frac{979}{18}\zae+\frac{(143+9\gamma)}{6\sqrt{3}}nW_1+\frac{(459+376\gamma)}{27\sqrt{3}}\zwe\notag\\&&+\gamma\left\{14+\frac{4\epsilon}{3}+25A_2-\frac{1507 }{18}\zae-\frac{(215+29\gamma)}{6\sqrt{3}}nW_1
-\frac{2(1174+169\gamma)}{27\sqrt{3}}\zwe\right\}\biggr]\\
T_2&=&\frac{3\sqrt{3}}{16}\biggl[14-\frac{16}{3}\epsilon+\frac{A_2}{3}-\frac{367}{18}\zae+\frac{115(1+\gamma)}{18\sqrt{3}}nW_1-\frac{(959-136\gamma)}{27\sqrt{3}}\zwe\\&&+\gamma\left\{\frac{32\epsilon}{3}+40A_2-\frac{382}{9}\zae+\frac{(511+53\gamma)}{6\sqrt{3}}nW_1-\frac{(2519-24\gamma)}{27\sqrt{3}}\zwe\right\}\biggr]\notag\\
 T_3&=&\frac{-9}{16}\biggl[\frac{8}{3}\epsilon+\frac{203A_2}{6}-\frac{625}{54}\zae-\frac{(105+15\gamma)}{18\sqrt{3}}nW_1-\frac{(403-114\gamma)}{81\sqrt{3}}\zwe\\&&+\gamma\left\{2-\frac{4\epsilon}{9}+\frac{55A_2}{2}-\frac{797}{54}\zae+\frac{(197+23\gamma)}{18\sqrt{3}}nW_1-\frac{(211-32\gamma)}{81\sqrt{3}}\zwe\right\}\biggr]\notag\\
 T_4&=&\frac{-9\sqrt{3}}{16}\biggl[2-\frac{8}{3}\epsilon+\frac{23A_2}{3}-44\zae-\frac{(37+\gamma)}{18\sqrt{3}}nW_1-\frac{(219+253\gamma)}{81\sqrt{3}}\zwe\\&&+\gamma\left\{4\epsilon+\frac{88}{27}\zae+\frac{(241+45\gamma)}{18\sqrt{3}}nW_1-\frac{(1558-126\gamma)}{81\sqrt{3}}\zwe\right\}\biggr]\notag\\ T_5&=&\frac{W_1}{2(a^2+b^2)^3}\biggl[(a\dot{x}+b\dot{y})\left\{3(ax+by)-(bx-ay)^2\right\}-2(x\dot{x}+y\dot{y})(ax+by)(a^2+b^2)\biggr]\label{eq:t5}\end{eqnarray}

The second order part $H_2$ of the corresponding Hamiltonian takes
the form
\begin{equation}
H_2=\frac{p_x^2+p_y^2}{2}+n(yp_x-xp_y)+Ex^2+F^2+Gxy
\end{equation}
To investigate the stability of the motion, as in Wittaker(1965), we
consider the following set of linear equations in the variables $x,
y$:
 \begin{equation}\label{eq:ax}
 \begin{array}{l c l}
 -\lambda P_x& = & \frac{\partial{H_2}}{\partial x}\\&&\\
 -\lambda P_y& = & \frac{\partial{H_2}}{\partial y}\\
 \text{i.e.}\quad AX&=&0
 \end{array} \quad
 \begin{array}{l c l }
 \lambda x& = & \frac{\partial{H_2}}{\partial P_x}\\&&\\
 \lambda y& = & \frac{\partial{H_2}}{\partial P_y}\\&&
 \end{array}
  \end{equation}
  \begin{equation}
  X=\left[\begin{array}{c}
  x\\
  y\\
  P_x\\
  P_y \end{array}\right] \quad \text{and}
 \quad
 A=\left[\begin{array}{c c c c}
  2E & G&\lambda& -n\\
G&2F&n&\lambda\\
  -\lambda& n& 1& 0\\
  -n & -\lambda& 0& 1\end{array}\right]
  \end{equation}

 Clearly $|A|=0$ implies that the characteristic equation
 corresponding to Hamiltonian $H_2$ is given by
 \begin{equation}
 \lambda^4+2(E+F+n^2)\lambda^2+4EF -G^2+n^4-2n^2(E+F)=0 \label{eq:ch}
 \end{equation}
 This is characteristic equation whose discriminant is
  \begin{equation}
 D=4(E+F+n^2)^2-4\bigl\{4EF-G^2+n^4-2n^2(E+F)\bigr\}
 \end{equation}
 Stability is assured  only when $D>0$.
 i.e
  \begin{eqnarray}
  \mu&<&\mu_{c_0}-0.221895916277307669\epsilon +2.1038871010983331 A_2
  \notag\\&+&
    0.493433373141671349\epsilon A_2 +0.704139054372097028 n W_1 +
    0.401154273957540929 n\epsilon W_1\notag
 \end{eqnarray}
 where $\mu_{c_0}=0.0385208965045513718$

 When $D>0$ the roots $\pm i \omega_1$ and $\pm i \omega_2$ ($\omega_1,\omega_2$ being the long/short -periodic frequencies) are related to each other as

 \begin{eqnarray}
   \omega_1^2+\omega_2^2&=& 1-\frac{\gamma \epsilon}{2}+\frac{3\gamma A_2}{2}+\frac{83\epsilon A_2}{12}+\frac{299\gamma\epsilon A_2}{144}-\frac{n W_1}{24\sqrt{3}}+\frac{5 \gamma n W_1}{8\sqrt{3}}-\frac{53 \epsilon n W_1}{54\sqrt{3}}\notag\\
  &&-\frac{5 \gamma^2 n W_1}{24\sqrt{3}}+\frac{173 \gamma \epsilon n W_1}{54\sqrt{3}}-\frac{3 \gamma^2 \epsilon n
  W_1}{36\sqrt{3}}\label{eq:w1+w2} \\
  \omega_1^2\omega_2^2&=&\frac{27}{16} -\frac{27\gamma^2}{16}+\frac{9\epsilon}{8}+\frac{9\gamma\epsilon}{8} -\frac{3\gamma^2\epsilon}{8}+\frac{117\gamma A_2}{16}-\frac{241\epsilon A_2}{32}+\frac{2515\gamma\epsilon A_2}{192}\notag\\
  &&+\frac{35n W_1}{16\sqrt{3}}-\frac{55 \sqrt{3}\gamma n W_1}{16}-\frac{5\sqrt{3} \gamma^2 n W_1}{4}-\frac{1277 \epsilon n W_1}{288\sqrt{3}}+\frac{5021 \gamma \epsilon n W_1}{288\sqrt{3}}+\frac{991 \gamma^2 \epsilon n W_1}{48\sqrt{3}}\label{eq:w1w2} \\
&&(0<\omega_2<\frac{1}{\sqrt{2}}<\omega_1<1)\notag\end{eqnarray}

From (~\ref{eq:w1+w2})and  (~\ref{eq:w1w2}) it may be noted that
$\omega_j (j=1,2)$ satisfy

\begin{eqnarray}
  \gamma^2&=& 1+\frac{4\epsilon}{9}-\frac{107\epsilon A_2}{27}+\frac{2\gamma \epsilon }{3}+\frac{1579\gamma\epsilon A_2}{324}-\frac{25nW_1}{27\sqrt{3}}-\frac{55\gamma nW_1}{9\sqrt{3}}+\frac{3809\epsilon nW_1}{486\sqrt{3}}+\frac{4961\gamma\epsilon nW_1}{486\sqrt{3}}\notag\\&&
  +\biggl(-\frac{16}{27}+\frac{32\epsilon}{243}+\frac{208 A_2}{81}-\frac{8\gamma A_2}{27}-\frac{4868\epsilon
  A_2}{729}-\frac{563\gamma\epsilon A_2}{243}\notag\\&&+\frac{296nW_1}{243\sqrt{3}}-\frac{10\gamma
  nW_1}{27\sqrt{3}}-\frac{15892\epsilon nW_1}{2187\sqrt{3}}-\frac{1864\gamma\epsilon
  nW_1}{729\sqrt{3}}\biggr)\omega_j^2\notag\\
  && +\biggl(\frac{16}{27}-\frac{32\epsilon}{243}-\frac{208 A_2}{81}-\frac{1880\epsilon
  A_2}{729}-\frac{2720nW_1}{2187\sqrt{3}}+\frac{49552\epsilon nW_1}{6561\sqrt{3}}-\frac{80\gamma\epsilon
  nW_1}{2187\sqrt{3}}\biggr)\omega_j^4
 \end{eqnarray}
 Alternatively, it can also be seen that if $u=\omega_1\omega_2$,
 then equation (~\ref{eq:w1w2}) gives
\begin{eqnarray}
  \gamma^2&=& 1+\frac{4\epsilon}{9}-\frac{107\epsilon A_2}{27}-\frac{25nW_1}{27\sqrt{3}}+\frac{3809\epsilon nW_1}{486\sqrt{3}}+\gamma\biggl(\frac{2\epsilon }{3}+\frac{1579\epsilon A_2}{324}-\frac{55\gamma nW_1}{9\sqrt{3}}+\frac{4961\gamma\epsilon
  nW_1}{486\sqrt{3}}\biggr)\notag\\&&
  +\biggl(-\frac{16}{27}+\frac{32\epsilon}{243}+\frac{208 A_2}{81}-\frac{1880\epsilon A_2}{729}+\frac{320nW_1}{243\sqrt{3}}-\frac{15856\epsilon nW_1}{2187\sqrt{3}}\biggr)u^2 \end{eqnarray}

Following the method for reducing $H_2$ to the normal form, as in
Whittaker(1965),use the transformation
\begin{equation}
 X=JT \quad  \text{where}  \quad X=\left[\begin{array}{c}
x\\y\\p_x\\p_y\end{array}\right],J=[J_{ij}]_{1\leq i\leq j \leq 4},\
T=\left[\begin{array}{c} Q_1\\Q_2\\p_1\\p_2\end{array}\right]
\end{equation}
\begin{equation}
P_i= (2 I_i\omega_i)^{1/2}\cos{\phi_i},  \quad Q_i= (\frac{2
I_i}{\omega_i})^{1/2}\sin{\phi_i}, \quad (i=1,2)
\end{equation}
The transformation changes the second order part of the Hamiltonian
into the normal form \begin{equation}
H_2=\omega_1I_1-\omega_2I_2\end{equation}

The general solution of the corresponding equations of motion are

\begin{equation}I_i=\text{const.}, \quad \phi_i=\pm \omega_i+\text{const},\  (i=1,2)\label{eq:I}\end{equation}
If the oscillations about $L_4$ are exactly linear, the
equation(~\ref{eq:I}) represent the integrals of motion and the
corresponding orbits will bi given by
\begin{eqnarray}x&=&J_{13}\sqrt{2\omega_1I_1}\cos{\phi_1}+J_{14}\sqrt{2\omega_2I_2}\cos{\phi_2}\label{eq:xb110}\\
y&=&J_{21}\sqrt{\frac{2I_1}{\omega_1}}\sin{\phi_1}+J_{22}\sqrt{\frac{2I_2}{\omega_2}}\sin{\phi_2}+J_{23}\sqrt{2I_1}{\omega_1}\cos{\phi_1}+J_{24}\sqrt{2I_2}{\omega_2}\sin{\phi_2}\label{eq:yb101}
\end{eqnarray}
where \begin{eqnarray}
J_{13}&=&\frac{l_1}{2\omega_1k_1}\left\{1-\frac{1}{2l_1^2}\left[\epsilon+\frac{45A_2}{2}-\frac{717A_2\epsilon}{36}+\frac{(67+19\gamma)}{12\sqrt{3}}nW_1
-\frac{(431-3\gamma)}{27\sqrt{3}}nW_1\epsilon\right]\right.\notag\\&&+\frac{\gamma}{2l_1^2}\left[3\epsilon-\frac{29A_2}{36}-\frac{(187+27\gamma)}{12\sqrt{3}}nW_1
-\frac{2(247+3\gamma)}{27\sqrt{3}}nW_1\epsilon\right]\notag\\&&-\frac{1}{2k_1^2}\left[\frac{\epsilon}{2}-3A_2-\frac{73A_2\epsilon}{24}+\frac{(1-9\gamma)}{24\sqrt{3}}nW_1
+\frac{(53-39\gamma)}{54\sqrt{3}}nW_1\epsilon\right]\notag\\&&-\frac{\gamma}{4k_1^2}\left[\epsilon-3A_2-\frac{299A_2\epsilon}{72}-\frac{(6-5\gamma)}{12\sqrt{3}}nW_1
-\frac{(266-93\gamma)}{54\sqrt{3}}nW_1\epsilon\right]\notag\\&&\left.+\frac{\epsilon}{4l_1^2k_1^2}\left[\frac{3A_2}{4}
+\frac{(33+14\gamma)}{12\sqrt{3}}nW_1\right]+\frac{\gamma\epsilon}{8l_1^2k_1^2}\left[\frac{347A_2}{36}
-\frac{(43-8\gamma)}{4\sqrt{3}}nW_1
\right]\right\}\end{eqnarray}
 \begin{eqnarray} J_{14}&=&\frac{l_2}{2\omega_2k_2}\left\{1-\frac{1}{2l_2^2}\left[\epsilon+\frac{45A_2}{2}-\frac{717A_2\epsilon}{36}+\frac{(67+19\gamma)}{12\sqrt{3}}nW_1
-\frac{(431-3\gamma)}{27\sqrt{3}}nW_1\epsilon\right]\right.\notag\\&&-\frac{\gamma}{2l_2^2}\left[3\epsilon-\frac{293A_2}{36}+\frac{(187+27\gamma)}{12\sqrt{3}}nW_1
-\frac{2(247+3\gamma)}{27\sqrt{3}}nW_1\epsilon\right]\notag\\&&-\frac{1}{2k_2^2}\left[\frac{\epsilon}{2}-3A_2-\frac{73A_2\epsilon}{24}+\frac{(1-9\gamma)}{24\sqrt{3}}nW_1
+\frac{(53-39\gamma)}{54\sqrt{3}}nW_1\epsilon\right]\notag\\&&+\frac{\gamma}{2k_2^2}\left[\epsilon-3A_2-\frac{299A_2\epsilon}{72}-\frac{(6-5\gamma)}{12\sqrt{3}}nW_1
-\frac{(268-9\gamma)}{54\sqrt{3}}nW_1\epsilon\right]\notag\\&&\left.-\frac{\epsilon}{4l_2^2k_2^2}\left[\frac{33A_2}{4}
+\frac{(1643-93\gamma)}{216\sqrt{3}}nW_1\right]+\frac{\gamma\epsilon}{4l_2^2k_2^2}\left[\frac{737A_2}{72}
-\frac{(13+2\gamma)}{\sqrt{3}}nW_1
\right]\right\}\end{eqnarray}
\begin{eqnarray} J_{21}&=&-\frac{4n\omega_1}{l_1k_1}\left\{1+\frac{1}{2l_1^2}\left[\epsilon+\frac{45A_2}{2}-\frac{717A_2\epsilon}{36}+\frac{(67+19\gamma)}{12\sqrt{3}}nW_1
-\frac{(413-3\gamma)}{27\sqrt{3}}nW_1\epsilon\right]\right.\notag\\&&-\frac{\gamma}{2l_1^2}\left[3\epsilon-\frac{293A_2}{36}+\frac{(187+27\gamma)}{12\sqrt{3}}nW_1
-\frac{2(247+3\gamma)}{27\sqrt{3}}nW_1\epsilon\right]\notag\\&&-\frac{1}{2k_1^2}\left[\frac{\epsilon}{2}-3A_2-\frac{73A_2\epsilon}{24}+\frac{(1-9\gamma)}{24\sqrt{3}}nW_1
+\frac{(53-39\gamma)}{54\sqrt{3}}nW_1\epsilon\right]\notag\\&&-\frac{\gamma}{4k_1^2}\left[\epsilon-3A_2-\frac{299A_2\epsilon}{72}-\frac{(6-5\gamma)}{12\sqrt{3}}nW_1
-\frac{(268-93\gamma)}{54\sqrt{3}}nW_1\epsilon\right]\notag\\&&\left.+\frac{\epsilon}{8l_1^2k_1^2}\left[\frac{33A_2}{4}+\frac{(68-10\gamma)}{24\sqrt{3}}nW_1\right]+\frac{\gamma\epsilon}{8l_1^2k_1^2}\left[\frac{242A_2}{9}
+\frac{(43-8\gamma)}{4\sqrt{3}}nW_1
\right]\right\}\end{eqnarray}
\begin{eqnarray} J_{22}&=&\frac{4n\omega_2}{l_2k_2}\left\{1+\frac{1}{2l_2^2}\left[\epsilon+\frac{45A_2}{2}-\frac{717A_2\epsilon}{36}+\frac{(67+19\gamma)}{12\sqrt{3}}nW_1
-\frac{(413-3\gamma)}{27\sqrt{3}}nW_1\epsilon\right]\right.\notag\\&&-\frac{\gamma}{2l_2^2}\left[3\epsilon-\frac{293A_2}{36}+\frac{(187+27\gamma)}{12\sqrt{3}}nW_1
-\frac{2(247+3\gamma)}{27\sqrt{3}}nW_1\epsilon\right]\notag\\&&+\frac{1}{2k_2^2}\left[\frac{\epsilon}{2}-3A_2-\frac{73A_2\epsilon}{24}+\frac{(1-9\gamma)}{24\sqrt{3}}nW_1
+\frac{(53-39\gamma)}{54\sqrt{3}}nW_1\epsilon\right]\notag\\&&-\frac{\gamma}{4k_2^2}\left[\epsilon-3A_2-\frac{299A_2\epsilon}{72}-\frac{(6-5\gamma)}{12\sqrt{3}}nW_1
-\frac{(268-93\gamma)}{54\sqrt{3}}nW_1\epsilon\right]\notag\\&&\left.+\frac{\epsilon}{4l_2^2k_2^2}\left[\frac{33A_2}{4}+\frac{(34+5\gamma)}{12\sqrt{3}}nW_1
\right]+\frac{\gamma\epsilon}{8l_2^2k_2^2}\left[\frac{75A_2}{2}+\frac{(43-8\gamma)}{4\sqrt{3}}nW_1
 \right]\right\}\end{eqnarray}
\begin{eqnarray} J_{23}&=&\frac{\sqrt{3}}{4\omega_1l_1k_1}\left\{2\epsilon+6A_2+\frac{37A_2\epsilon}{2}-\frac{(13+\gamma)}{2\sqrt{3}}nW_1
+\frac{2(79-7\gamma)}{9\sqrt{3}}nW_1\epsilon\right.\notag\\&&-\gamma\left[6+\frac{2\epsilon}{3}+13A_2-\frac{33A_2\epsilon}{2}+\frac{(11-\gamma)}{2\sqrt{3}}nW_1
-\frac{(186-\gamma)}{9\sqrt{3}}nW_1\epsilon\right]\notag\\&&+\frac{1}{2l_1^2}\left[51A_2+\frac{(14+8\gamma)}{3\sqrt{3}}nW_1\right]-\frac{\epsilon}{k_1^2}\left[3A_2
+\frac{(19+6\gamma)}{6\sqrt{3}}nW_1\right]\notag\\&&-\frac{\gamma}{2l_1^2}\left[6\epsilon+135A_2-\frac{808A_2\epsilon}{9}-\frac{(67+19\gamma)}{2\sqrt{3}}nW_1
-\frac{(755+19\gamma)}{9\sqrt{3}}nW_1\epsilon\right]\notag\\&&-\frac{\gamma}{2k_1^2}\left[3\epsilon-18A_2-\frac{55A_2\epsilon}{4}-\frac{(1-9\gamma)}{4\sqrt{3}}nW_1
+\frac{(923-60\gamma)}{12\sqrt{3}}nW_1\epsilon\right]\notag\\&&\left.+\frac{\gamma\epsilon}{8l_1^2k_1^2}\left[\frac{9A_2}{2}
+\frac{(34-5\gamma)}{2\sqrt{3}}nW_1\right]\right\}\qquad\end{eqnarray}
\begin{eqnarray}                 J_{24}&=&\frac{\sqrt{3}}{4\omega_2l_2k_2}\left\{2\epsilon+6A_2+\frac{37A_2\epsilon}{2}-\frac{(13+\gamma)}{2\sqrt{3}}nW_1
+\frac{2(79-7\gamma)}{9\sqrt{3}}nW_1\epsilon\right.\label{eq:j24}\notag\\
&&-\gamma\left[6+\frac{2\epsilon}{3}+13A_2-\frac{33A_2\epsilon}{2}+\frac{(11-\gamma)}{2\sqrt{3}}nW_1
-\frac{(186-\gamma)}{9\sqrt{3}}nW_1\epsilon\right]\notag\\&&-\frac{1}{2l_2^2}\left[51A_2+\frac{(14+8\gamma)}{3\sqrt{3}}nW_1\right]-\frac{\epsilon}{k_2^2}\left[3A_2
+\frac{(19+6\gamma)}{6\sqrt{3}}nW_1\right]\notag\\&&-\frac{\gamma}{2l_2^2}\left[6\epsilon+135A_2-\frac{808A_2\epsilon}{9}-\frac{(67+19\gamma)}{2\sqrt{3}}nW_1
-\frac{(755+19\gamma)}{9\sqrt{3}}nW_1\epsilon\right]\notag\\&&-\frac{\gamma}{2k_1^2}\left[3\epsilon-18A_2-\frac{55A_2\epsilon}{4}-\frac{(1-9\gamma)}{4\sqrt{3}}nW_1
+\frac{(923-60\gamma)}{12\sqrt{3}}nW_1\epsilon\right]\notag\\&&\left.-\frac{\gamma\epsilon}{4l_1^2k_1^2}\left[\frac{99A_2}{2}
+\frac{(34-5\gamma)}{2\sqrt{3}}nW_1\right]\right\}\end{eqnarray}
with $l_j^2=4\omega_j^2+9,(j=1,2)$ and $ k_1^2=2\omega_1^2-1,
k_2^2=-2\omega_2^2+1 $.
\section{Conclusion}
Using Wittaker(1965) method we find that the second order part
 $H_2$ of the Hamiltonian is transformed into the normal
form $H_2=\omega_1I_1-\omega_2I_2.$

\thanks{\bf Acknowledgment:}{ We are thankful to D.S.T. Government of India, New Delhi for sanctioning a project DST/MS/140/2K dated 02/01/2004 on this topic. We are also thankful to IUCAA Pune for providing  financial assistance for visiting library and  computer facility.}

\end{document}